\documentclass[12pt]{article}
\usepackage{amsfonts}
\usepackage{latexsym, amssymb, amsmath, amscd, amsfonts,mathrsfs}

\newtheorem{thm}{Theorem}[section]

\numberwithin{equation}{section}

\numberwithin{equation}{section}

\def\qed{\hfill \rule{4pt}{7pt}}

\def\pf{\noindent {\it Proof.} }

\begin{document}
\begin{center}
{\large {\bf A Franklin Type Involution for Squares}}

\vskip 6mm {\small William Y.C. Chen$^1$ and  Eric H. Liu$^2$ \\[%
2mm] Center for Combinatorics, LPMC-TJKLC\\
Nankai University, Tianjin 300071,
P.R. China \\[3mm]
$^1$chen@nankai.edu.cn,  $^2$eric@cfc.nankai.edu.cn \\[0pt%
] }
\end{center}

\vskip 6mm \noindent{\bf Abstract.} We find an involution
as a combinatorial proof of a Ramanujan's partial
theta identity. Based on this involution, we obtain
a Franklin type involution for
squares in the sense that the classical Franklin involution
provides a combinatorial interpretation of
Euler's pentagonal number theorem. This Franklin type involution
 can be considered as a solution to a problem
 proposed by Pak concerning the parity of the number of partitions of
$n$ into distinct parts with the smallest part being odd.
Using a weighted form of our
involution, we give a combinatorial proof of a weighted
partition theorem derived by Alladi from
 Ramanujan's partial theta identity. This answers a question
of Berndt, Kim and Yee. Furthermore, through a different
weight assignment, we find combinatorial interpretations for
another partition theorem
derived by Alladi from
a partial theta identity of Andrews.  Moreover, we
obtain a partition theorem based on Andrews' identity and provide
a combinatorial proof by certain weight assignment for
our involution. A specialization of our partition theorem  is relate to an identity of Andrews
 concerning partitions into distinct nonnegative parts with the smallest part being even.
 Finally, we give a more general form of our partition theorem which in return corresponds
 to a generalization of Andrews' identity.

\noindent {\bf Keywords}: Franklin type involution, Ramanujan's partial theta identity, Andrews' identity.

\noindent {\bf AMS Mathematical Subject Classifications}: 05A17,
11P81

\section{ Introduction}

The main result of this paper is a Franklin type involution for squares which is related to Ramanujan's partial theta identity and an identity of Andrews.
As applications of this involution, we answer a question of Pak \cite[p. 31,5.2.7]{pak06} on the parity
of the number of partitions of $n$ into distinct parts with the smallest
part being odd, and we give a solution to a problem proposed by Berndt,
Kim and Yee \cite{berndt09}
by providing a combinatorial interpretation of a partition theorem derived by
Alladi \cite{alld091} from Ramanujan's partial theta identity. Furthermore,
we obtain a partition theorem based on Andrews' identity. A specialization
of this theorem is related to an identity of Andrews on partitions into
distinct nonnegative parts with the smallest part being even.  Finally, we
find a more general form of our partition theorem which in return corresponds to a
generalization of Andrews' identity.

Recall that the celebrated involution of Franklin gives a
combinatorial interpretation of
Euler's pentagonal number theorem as stated below
\begin{equation}\label{euler}
(q;q)_{\infty}=1+\sum\limits_{k=1}^{\infty}(-1)^k(q^{k(3k-1)/2}+q^{k(3k+1)/2}),
\end{equation}
where the $q$-shifted factorial is defined by
$$(a;q)_n:=(1-a)(1-aq)\cdots(1-aq^{n-1}), \quad n\geq 1,$$
and
$$(a;q)_{\infty}=\lim\limits_{n\rightarrow \infty}(a;q)_n, \quad |q|<1.$$
Let $D$ denote the set of integer partitions
into distinct parts, and let $D(n)$ denote the set of
partitions of $n$ into distinct parts. The relation  (\ref{euler}) can be
 reinterpreted as the following number-theoretic identity
\begin{equation}\label{euler2}
\sum\limits_{ \lambda \in D(n)}(-1)^{\ell(\lambda)}=
\left\{
\begin{array}{ll}
(-1)^k,  & \mbox{if}\ n=k(3k\pm 1)/2,\\[6pt]
0, & \mbox{otherwise},
\end{array}
\right.
\end{equation}
where
$\ell(\lambda)$ denotes the
number of parts of  $\lambda$.

The Franklin type involution for squares will be concerned with
the set of partitions of a non-negative integer
into distinct parts with the smallest part being odd.
Let us use $P_{do}$ to denote the set of such partitions, and use $P_{do}(n)$
to denote the set of such partitions of $n$.
To be more specific, we obtain the following number-theoretic identity which is
analogous to (\ref{euler2}),
\begin{equation}\label{euler3}
\sum\limits_{\lambda \in P_{do}(n)}(-1)^{\ell(\lambda)}=\left\{
\begin{array}{ll}
(-1)^k,  & \mbox{if}\ n=k^2,\\[6pt]
0, & \mbox{otherwise}.
\end{array}
\right.
\end{equation}
It is clear that (\ref{euler3}) implies a solution of a problem
of Pak \cite[p. 31,5.2.7]{pak06} concerning the parity of the number of
partitions in $P_{do}(n)$. More precisely, he posed the problem
 of finding an explicit involution to justify  that
the number of elements in $P_{do}(n)$ is odd if and
only if $n$ is a square.

Moreover, for various
weight assignments $\omega(\lambda)$  to partitions $\lambda \in P_{do}$,
 our involution turns out to be sign-revering and
weight-preserving. This property leads us to several
number-theoretic identities
of the following form:
\begin{equation}\label{euler4}
\sum\limits_{\lambda \in P_{do}(n)}\omega(\lambda)=\left\{
\begin{array}{ll}
(-a)^k,  & \mbox{if}\ n=k^2,\\[6pt]
0, & \mbox{otherwise}.
\end{array}
\right.
\end{equation}
The first case is related to a problem proposed by
Berndt, Kim and Yee \cite{berndt09} concerning
a combinatorial interpretation of a number-theoretic identity
derived by Alladi \cite{alld091} from the following
Ramanujan's partial theta identity from Ramanujan's lost notebook
\cite[p. 38]{raman88},
\begin{equation}\label{iden1}
1+\sum\limits_{k=1}^{\infty}\frac{(-q;q)_{k-1}(-a)^kq^{k(k+1)/2}}
{(aq^2;q^2)_k}=\sum\limits_{k=0}^{\infty}(-a)^kq^{k^2}.
\end{equation}
By giving a weight function in terms of the gaps between the parts
 of partitions in $P_{do}$, Alladi \cite{alld091} derived a partition theorem
 in the above form, see, Section 4.
Though Berndt, Kim and Yee \cite{berndt09} have found a bijective proof Ramanujan's identity (\ref{iden1}),
it is not clear whether their involution can imply a combinatorial interpretation
of Alladi's weighted partition theorem. As will be seen, our Franklin type involution indeed gives a combinatorial proof
of Alladi's partition theorem.

The second case is concerned with a weighted partition theorem
obtained by Alladi \cite{alld092} from the following partial theta identity
of Andrews \cite[p. 157]{andr76},
\begin{equation}\label{iden2}
\sum\limits_{n=0}^{\infty}q^{2n}(q^{2n+2};q^2)_{\infty}(aq^{2n+1};q^2)_{\infty}
=\sum\limits_{k=0}^{\infty}(-a)^kq^{k^2}.
\end{equation}
By giving a weight function in terms of the odd parts
of partitions in $P_{do}$, Alladi \cite{alld092} derived a partition theorem
in the form of (\ref{euler4}), see, Section 5.
It turns out that our involution also applies to this partition
theorem with respect to a different weight assignment.

As the third application of our involution, we give a combinatorial
 proof of a number-theoretic
 theorem on partitions
into distinct parts with smallest part being even derived from
Andrews' identity (\ref{iden2}). Moreover, we note that
a special case of this partition theorem is related to an identity of
Andrews, first proposed as a problem in  \cite{andr72}. A generating function proof was given by
Stenger \cite{sten73}.

To conclude this paper, we extend our involution  to derive
 a more general  identity
\begin{equation}\label{iden3}
\begin{split}
&\sum\limits_{n=0}^{\infty}q^{2mn}(q^{2mn+2m};q^{2m})_{\infty}(aq^{2mn+1};q^2)_{\infty}\\
&=
1+\sum\limits_{k=1}^{\infty}(-a)^kq^{k^2}\prod\limits_{j=1}^k(1+q^{2j}+q^{4j}
+\cdots+q^{2(m-1)j}),
\end{split}
\end{equation}
which reduces to the following identity of Andrews
 \cite[p. 157]{andr76} when setting $a=-1$,
\begin{equation}\label{add7}
\begin{split}
&\sum\limits_{n=0}^{\infty}q^{2mn}(q^{2mn+2m};q^{2m})_{\infty}(-q^{2mn+1};q^2)_{\infty}\\
&=1+\sum\limits_{k=1}^{\infty}q^{k^2}\prod\limits_{j=1}^k(1+q^{2j}
+q^{4j}+\cdots+q^{2(m-1)j}).
\end{split}
\end{equation}
Notice that (\ref{add7}) is a generalization of  (\ref{iden2}).

\section{An involution for Ramanujan's identity}

In this section,
we shall construct an  involution which leads to
a combinatorial proof of Ramanujan's partial theta identity (\ref{iden1}).
This involution serves as a
crucial step in the Franklin type involution presented in the next section
which can be viewed as a
bijective proof of Alladi's partition theorem derived from (\ref{iden1}) with respect to certain weight
assignment.

Let $D_k$ be the set of partitions $\pi$ into $k$
distinct parts with the smallest part being $1$ such that
$\pi_i-\pi_{i+1}\leq 2$, and
let $E_k$ denote the set of partitions
$\sigma$ with even parts not exceeding $2k$,
that is, each $\sigma_i$ is even and $\sigma_1 \leq 2k$.
We are going to establish an involution on $D_k\times E_k$.
Throughout this paper, $T_k$ standards for the triangular partition
$(2k-1,2k-3,\ldots,3,1)$.

\begin{thm}
There exists an involution on the set $D_k\times E_k$ under which the pair of
partitions $(T_k,\emptyset)$ remains invariant.
\end{thm}

To construct the desired involution on $D_k\times E_k$, we introduce a
statistic called the modular leg hook length of a partition in $D_k$.
Adopting the notation  in \cite{pak06},
we use $[\lambda]_2$ to denote the 2-modular diagram  of a partition $\lambda$
defined to be a Young diagram filled with  $1$ or $2$ such that
the last cell of row $i$ is filled with $1$ if and only if $\lambda_i$ is odd.
Given a partition
$\pi=(\pi_1,\pi_2,\ldots,\pi_k) \in D_k$, let us consider the $2$-modular diagram.
Suppose that $\pi_i$ is an even part other than the largest part,
we can associate it with a  modular leg hook  $H_i$ which
consists of the squares in the $i$-th row in the $2$-modular diagram and the
squares in first column above the $i$-th row. For a modular leg hook
$H_i$, the  length of this hook, denoted by $|H_i|$, is defined to be the sum of the
numbers filled in the hook, and the  height of it is the number of squares
in the first column of this hook.

We are now ready to describe the construction of the
involution on $D_k\times E_k$. Let us denote this involution by $\varphi$.

{\noindent \bf The involution $\varphi$ on $D_k\times E_k$:} Given a
pair of partitions $(\pi, \sigma)\in D_k\times E_k$, represent
$\pi$ and $\sigma$ by their $2$-modular diagrams, respectively.
 In fact, the desired involution consists of two involutions.

{\noindent \it Part A:}
We have the following two cases.
\begin{itemize}
\item[(1)]
Suppose that there exists a modular leg hook in $\pi$ such
that after the deletion of this hook the resulting partition is in $D_{k-1}$, then we choose
such a hook  with  maximum height and denote it by $H_r(\pi)$. If $|H_r(\pi)|\geq \sigma_1$. Then delete
$H_r(\pi)$ from
$\pi$ and add it to  $\sigma$ as a new part.  We denote the
resulting partitions by $\pi^{*}$ and $\sigma^{*}$, respectively. Since $|H_r(\pi)| \leq 2k-2$,
we have $(\pi^{*},\sigma^{*})
\in D_{k-1}\times E_{k-1}$.

\item[(2)] Suppose that either there is the modular leg hook $H_r(\pi)$ in $\pi$ and $|H_r(\pi)|< \sigma_1$
or there does not exist the modular leg hook $H_r(\pi)$ in $\pi$ and $\pi_1+2< \sigma_1$. Then insert $\sigma_1$ into $\pi$ as a modular leg hook in the $2$-modular diagram of $\pi$. To be precise, this operation
can be described as follows.
 Let $i$ be the largest positive integer such that
$\sigma_1-2i>\pi_{i+1}$, that is, for $j>i$ we have
$\sigma_1-2j<\pi_{j+1}$. Then we add $2$ to each of  the
first $i$ parts $\pi_1, \pi_2, \ldots, \pi_i$, and insert
$\sigma_1-2i$ as a new part before the part $\pi_{i+1}$.
Since
$\sigma_1\leq 2k-2$ and any two consecutive parts of $\pi$ differ by at most $2$,
the resulting pair of partitions, denoted by $(\pi^{*},\sigma^{*})$,
belongs to  $D_{k+1}\times E_k$. Furthermore, there exists the modular leg hook $H_r(\pi^*)$ in $\pi^*$
and $|H_r(\pi^*)|\geq \sigma_1^*$.
\end{itemize}
Below is an example.
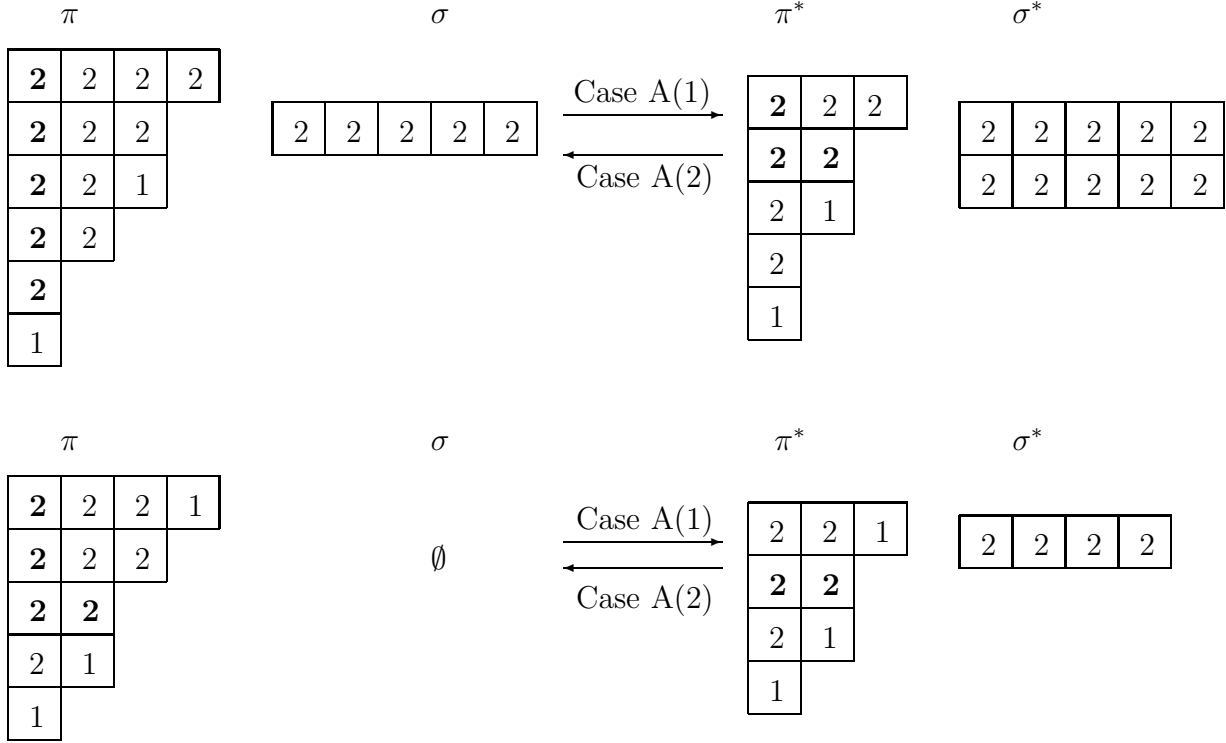
\begin{center}
\begin{picture}(460,130)
\put(0,120){\line(1,0){80}}\put(0,100){\line(1,0){80}}
\put(0,80){\line(1,0){60}}\put(0,60){\line(1,0){60}}
\put(0,40){\line(1,0){40}}\put(0,20){\line(1,0){20}}
\put(0,0){\line(1,0){20}}\put(0,120){\line(0,-1){120}}
\put(20,120){\line(0,-1){120}}\put(40,120){\line(0,-1){80}}
\put(60,120){\line(0,-1){60}}\put(80,120){\line(0,-1){20}}
\multiput(28,105)(20,0){3}{2}\multiput(28,85)(20,0){2}{2}
\multiput(28,65)(20,0){1}{2}\put(28,45){2}\multiput(8,105)(0,-20){5}{\bf 2}
\put(8,5){1}\put(48,65){1}

\put(100,100){\line(1,0){100}}\put(100,80){\line(1,0){100}}
\multiput(100,100)(20,0){6}{\line(0,-1){20}}
\multiput(108,85)(20,0){5}{2}

\put(210,95){\vector(1,0){60}}\put(214,100){Case A(1)}
\put(270,80){\vector(-1,0){60}}\put(215,68){Case A(2)}

\put(280,110){\line(1,0){60}}\put(280,90){\line(1,0){60}}
\put(280,70){\line(1,0){40}}\put(280,50){\line(1,0){40}}
\put(280,30){\line(1,0){20}}\put(280,10){\line(1,0){20}}
\put(280,110){\line(0,-1){100}}\put(300,110){\line(0,-1){100}}
\put(320,110){\line(0,-1){60}}\put(340,110){\line(0,-1){20}}
\put(288,95){\bf 2}\put(308,95){2}\put(325,95){2}
\put(288,75){\bf 2}\put(308,75){\bf 2}
\put(288,55){2}\put(308,55){1}\put(288,35){2}\put(288,15){1}

\put(360,100){\line(1,0){100}}\put(360,80){\line(1,0){100}}
\put(360,60){\line(1,0){100}}\multiput(360,100)(20,0){6}{\line(0,-1){40}}
\multiput(368,85)(20,0){5}{2}\multiput(368,65)(20,0){5}{2}

\put(20,130){$\pi$}\put(160,130){$\sigma$}\put(290,130){$\pi^{*}$}\put(380,130){$\sigma^{*}$}

\end{picture}
\end{center}

\begin{center}
\begin{picture}(460,110)

\put(0,100){\line(1,0){80}}\put(0,80){\line(1,0){80}}
\put(0,60){\line(1,0){60}}\put(0,40){\line(1,0){40}}
\put(0,20){\line(1,0){40}}\put(0,0){\line(1,0){20}}
\put(0,100){\line(0,-1){100}}
\put(20,100){\line(0,-1){100}}\put(40,100){\line(0,-1){80}}
\put(60,100){\line(0,-1){40}}\put(80,100){\line(0,-1){20}}
\multiput(28,85)(20,0){2}{2}\multiput(28,65)(20,0){2}{2}
\put(28,45){{\bf 2}}\multiput(8,85)(0,-20){3}{{\bf 2}}
\put(8,5){1}\put(28,25){1}\put(68,85){1}\put(8,25){2}

\put(160,65){$\emptyset$}

\put(210,75){\vector(1,0){60}}\put(215,80){Case A(1)}
\put(270,65){\vector(-1,0){60}}\put(215,50){Case A(2)}

\put(280,90){\line(1,0){60}}\put(280,70){\line(1,0){60}}
\put(280,50){\line(1,0){40}}\put(280,30){\line(1,0){40}}
\put(280,10){\line(1,0){20}}
\put(280,90){\line(0,-1){80}}
\put(300,90){\line(0,-1){80}}\put(320,90){\line(0,-1){60}}
\put(340,90){\line(0,-1){20}}
\multiput(288,75)(20,0){2}{2}
\put(288,35){2}\multiput(288,55)(20,0){2}{{\bf 2}}
\put(328,75){1}\put(308,35){1}\put(288,15){1}

\put(360,85){\line(1,0){80}}\put(360,65){\line(1,0){80}}
\multiput(360,85)(20,0){5}{\line(0,-1){20}}
\multiput(368,70)(20,0){4}{2}

\put(20,110){$\pi$}\put(160,110){$\sigma$}\put(290,110){$\pi^{*}$}\put(380,110){$\sigma^{*}$}

\end{picture}
\end{center}

{\noindent \it Part B:} Suppose that there does not exist the modular leg hook $H_r(\pi)$ in $\pi$ and
$\sigma_1\leq \pi_1+2$. If $\pi$ has even parts,
then we choose the largest even part of $\pi$, and  denote it by $\pi_r$.
We consider the following two cases.

\begin{itemize}
\item[(1)]If $\pi_r\geq \sigma_1$,
then remove the part $\pi_r$ in $\pi$  and add it to $\sigma$.
We denote the resulting partitions by $\pi^{*}$ and $\sigma^{*}$, respectively.
Since $\pi_r \leq 2k-2$, we see that  $(\pi^{*},\sigma^{*})\in
D_{k-1}\times E_{k-1}$.
\item[(2)]If $\pi_r< \sigma_1$,
then remove  the part $\sigma_1$ in $\sigma$ and add it to  $\pi$. Denote the
resulting partitions by $\pi^{*}$ and $\sigma^{*}$, respectively. Since $\pi_i-\pi_{i+1}\leq 2$ for each
$i$,  $\sigma_1$ can be inserted either between two odd parts of $\pi$
or at the top of $\pi$. Therefore, $(\pi^{*},\sigma^{*})\in D_{k+1}\times E_{k}$.
\end{itemize}
Here is an example.
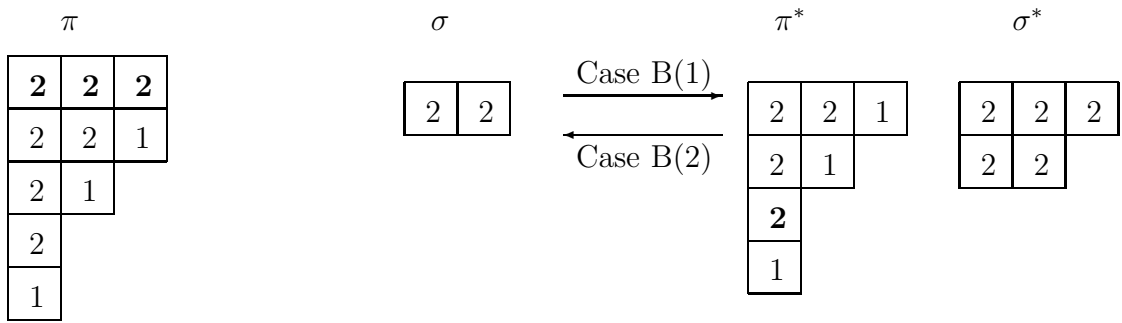
\begin{center}
\begin{picture}(460,110)
\put(0,100){\line(1,0){60}}\put(0,80){\line(1,0){60}}
\put(0,60){\line(1,0){60}}\put(0,40){\line(1,0){40}}
\put(0,20){\line(1,0){20}}\put(0,0){\line(1,0){20}}
\put(0,100){\line(0,-1){100}}
\put(20,100){\line(0,-1){100}}\put(40,100){\line(0,-1){60}}
\put(60,100){\line(0,-1){40}}
\multiput(8,65)(20,0){2}{2}
\put(48,65){1}\multiput(8,85)(20,0){3}{{\bf 2}}
\put(8,45){2}\put(28,45){1}\put(8,25){2}\put(8,5){1}

\put(150,90){\line(1,0){40}}\put(150,70){\line(1,0){40}}
\multiput(150,90)(20,0){3}{\line(0,-1){20}}\multiput(158,75)(20,0){2}{2}

\put(210,85){\vector(1,0){60}}\put(215,90){Case B(1)}
\put(270,70){\vector(-1,0){60}}\put(215,58){Case B(2)}

\put(280,90){\line(1,0){60}}\put(280,70){\line(1,0){60}}
\put(280,50){\line(1,0){40}}
\put(280,30){\line(1,0){20}}\put(280,10){\line(1,0){20}}
\put(280,90){\line(0,-1){80}}
\put(300,90){\line(0,-1){80}}
\put(320,90){\line(0,-1){40}}\put(340,90){\line(0,-1){20}}
\multiput(288,75)(20,0){2}{2}
\put(328,75){1}
\put(288,55){2}\put(308,55){1}\put(288,35){{\bf 2}}\put(288,15){1}

\put(360,90){\line(1,0){60}}\put(360,70){\line(1,0){60}}
\put(360,50){\line(1,0){40}}\multiput(360,90)(20,0){3}{\line(0,-1){40}}
\put(420,90){\line(0,-1){20}}\multiput(368,75)(20,0){3}{2}\multiput(368,55)(20,0){2}{2}

\put(20,110){$\pi$}\put(160,110){$\sigma$}\put(290,110){$\pi^{*}$}\put(380,110){$\sigma^{*}$}

\end{picture}
\end{center}

Finally, we are left with the case when
 $\pi$ has no even parts and $\sigma$ is the empty partition. In this situation, there is
 only one pair of partitions $(T_k, \emptyset)$, which is defined as the fixed point of the involution.

 It is straightforward to check that the above correspondence is
  an involution. Except for the fixed point, the mapping
changes the number of
even parts of $\pi$ by $1$ and preserves the number of odd parts at the same time.
Indeed, the above involution serves as a combinatorial proof of Ramanujan's partial
theta identity (\ref{iden1}).

{\noindent \it Proof of (\ref{iden1})}: Note that the generating function for
partitions $\pi \in D_k$ equals
\begin{equation}  \label{qq}
(-q;q)_{k-1}q^{k(k+1)/2},
\end{equation}
and the generating
function for partitions $\sigma \in E_k$ equals
$$\frac{1}{(q^2;q^2)_{k}}.$$
Thus the left hand side of (\ref{iden1}) corresponds to
partitions $(\pi, \sigma)\in \bigcup \limits_{k=0}^{\infty}D_k\times E_k$
with weight $(-1)^{\ell(\pi)}a^{\ell(\pi)+\ell(\sigma)}$.
Notice that the involution $\varphi$ changes the parity of
$\ell(\pi)$ and preserves the quantity $\ell(\pi)+\ell(\sigma)$.
The fixed point $(T_k, \emptyset)$
corresponds to  the right hand side of (\ref{iden1}).
In view of the involution $\varphi$, we obtain the identity
(\ref{iden1}).  \qed

In comparison, we note  that the bijective proof of (\ref{iden1}) given
by Berndt, Kim and Yee \cite{berndt09} is based on the interpretation of the
numerator (\ref{qq})
 in terms of parity sequences.

\section{A Franklin type involution for squares}

In this section, we shall construct
a Franklin type involution for squares where the involution $\varphi$ on $D_k\times E_k$ given in
the previous section serves as the main ingredient.
This involution is based on $P_{do}(n)$, namely, the set of partitions of
$n$ into distinct parts with the smallest part being odd, and it gives an
answer to a problem proposed by Pak on the characterization of the parity of
$|P_{do}(n)|$.

For partitions
$\lambda=(\lambda_1,\lambda_2,\ldots)$ and
$\mu=(\mu_1,\mu_2,\ldots),$ define $\lambda+\mu$ to be
the partition $(\lambda_1+\mu_1,\lambda_2+\mu_2,\ldots)$.
Denote the number of even (resp. odd) parts of $\lambda$ by
$\ell_e(\lambda)$ (resp. $\ell_o(\lambda)$).
Our involution  on $\lambda \in P_{do}(n)$, denoted by $\Psi$,
can be stated as follows.

{\noindent  Step $1$. Extraction of parts from $\lambda$:}
For a partition $\lambda=(\lambda_1,\lambda_2,\ldots,\lambda_k)\in
P_{do}(n)$, represent it by the $2$-modular diagram
$[\lambda]_2$, from which we can construct a pair of
partitions $(\pi, \sigma)\in D_k\times E_k$.
Initially, set $\pi=\lambda$, $\sigma=\emptyset$ and $t=k$.
 Then iterate the following procedure until $t=1$:
\begin{itemize}
\item[$\bullet$] Suppose  that there exists $i$ such that
$\pi_t-\pi_{t+1}=2i+r_t$, where $1\leq  r_t\leq 2$ and $\pi_{k+1}$ is defined to be $0$.
\item[$\bullet$] Subtract $2i$ from each of the parts $\pi_1, \pi_2,
\ldots, \pi_t$;
\item[$\bullet$] Rearrange the parts to form a new partition $\pi$
and add $i$ parts of size $2t$ to $\sigma$. Replace $t$ by $t-1$.
\end{itemize}
When $t=1$, we get a pair of partitions $(\pi, \sigma)\in D_k\times E_k$.
It is clear that
$$\ell(\lambda)=\ell(\pi),\quad
\ell_e(\lambda)=\ell_e(\pi), \quad \ell_o(\lambda)=\ell_o(\pi).$$
Here is an example.

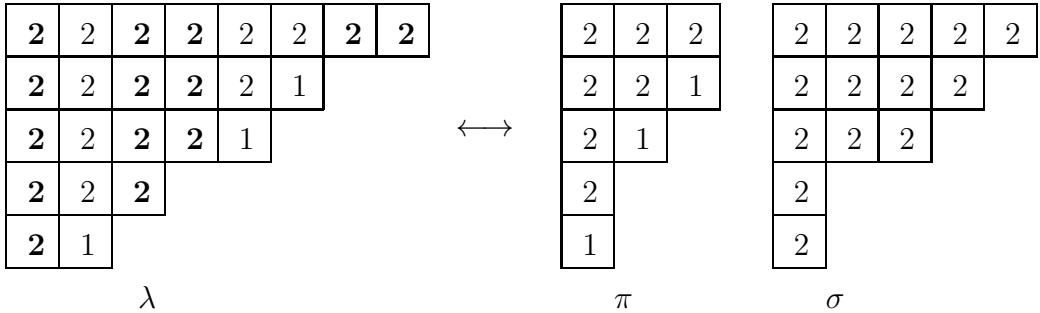
\begin{center}
\begin{picture}(390,130)
\put(0,130){\line(1,0){160}}\put(0,110){\line(1,0){160}}
\put(0,90){\line(1,0){120}}\put(0,70){\line(1,0){100}}
\put(0,50){\line(1,0){60}}\put(0,30){\line(1,0){40}}
\put(0,130){\line(0,-1){100}}\put(20,130){\line(0,-1){100}}
\put(40,130){\line(0,-1){100}}\put(60,130){\line(0,-1){80}}
\put(80,130){\line(0,-1){60}}\put(100,130){\line(0,-1){60}}
\put(120,130){\line(0,-1){40}}\put(140,130){\line(0,-1){20}}
\put(160,130){\line(0,-1){20}}\multiput(8,115)(0,-20){5}{\bf 2}
\multiput(28,115)(0,-20){4}{2}\multiput(48,115)(0,-20){4}{\bf 2}
\multiput(68,115)(0,-20){3}{\bf 2}\multiput(88,115)(0,-20){2}{2}
\multiput(108,115)(0,-20){1}{2}\put(128,115){\bf 2}
\put(148,115){\bf 2}\put(108,95){1}\put(88,75){1}
\put(28,35){1}

\put(170,80){$\longleftrightarrow$}

\put(210,130){\line(1,0){60}}\put(210,110){\line(1,0){60}}
\put(210,90){\line(1,0){60}}\put(210,70){\line(1,0){40}}
\put(210,50){\line(1,0){20}}\put(210,30){\line(1,0){20}}
\put(210,130){\line(0,-1){100}}\put(230,130){\line(0,-1){100}}
\put(250,130){\line(0,-1){60}}\put(270,130){\line(0,-1){40}}
\put(160,130){\line(0,-1){20}}\multiput(218,115)(0,-20){4}{2}
\multiput(238,115)(0,-20){2}{2}
\put(258,115){2}\put(258,95){1}
\put(238,75){1}\put(218,35){1}

\put(290,130){\line(1,0){100}}\put(290,110){\line(1,0){100}}
\put(290,90){\line(1,0){80}}\put(290,70){\line(1,0){60}}
\put(290,50){\line(1,0){20}}\put(290,30){\line(1,0){20}}
\put(290,130){\line(0,-1){100}}\put(310,130){\line(0,-1){100}}
\put(330,130){\line(0,-1){60}}\put(350,130){\line(0,-1){60}}
\put(370,130){\line(0,-1){40}}\put(390,130){\line(0,-1){20}}
\multiput(298,115)(0,-20){5}{2}
\multiput(318,115)(0,-20){3}{2}\multiput(338,115)(0,-20){3}{2}
\multiput(358,115)(0,-20){2}{2}
\put(378,115){2}

\put(50,15){$\lambda$}\put(230,15){$\pi$}\put(310,15){$\sigma$}

\end{picture}
\end{center}

{\noindent  Step 2. Apply the involution $\varphi$ on $D_k\times E_k$:}
For a pair of partitions $(\pi, \sigma)\in D_k\times E_k$,
use the involution $\varphi$ to generate a pair of partitions
$(\pi^{*},\sigma^{*})
\in D_{k+1}\times E_k$ or $
D_{k-1}\times E_{k-1}$.

{\noindent  Step 3. Insertion of parts of $\sigma^{*}$ to $\pi^{*}$:}
For a pair of partitions $(\pi^{*},\sigma^{*})
\in D_{k+1}\times E_k$ or $
D_{k-1}\times E_{k-1}$, consider their $2$-modular
diagrams.  Let $\lambda^{*}=\pi^*+c_2(\sigma^*)$.
Clearly, $\lambda^* \in P_{do}(n)$,
where $c_{2}(\sigma^*)$ denotes the $2$-modular conjugate partition
obtained  from
$[\sigma^*]_{2}$. It is easily seen that
$$\ell(\lambda^*)=\ell(\pi^*), \quad
\ell_e(\lambda^*)=\ell_e(\pi^*), \quad \ell_o(\lambda^*)=\ell_o(\pi^*).$$

Based on the above procedure, we can see that
the mapping $\Psi$ is  a bijection. Moreover, it possesses
the following property
\begin{equation}\label{add8}
\ell(\lambda^*)=\ell(\lambda)\pm 1,\quad
\ell_e(\lambda^*)=\ell_e(\lambda)\pm 1,\quad \ell_o(\lambda^*)=\ell_o(\lambda),
\end{equation}
where the $\pm$ sign means either plus or minus.
In other words, $\Psi$ changes the parity of the number of parts.
It is easy to check that only when
$n$ is a square,
say, $n=k^2$, there is exactly one partition which is undefined for $\Psi$,
that is, $\lambda=(2k-1,2k-3,\ldots,3,1)$. Hence we obtain an identity in
the spirit of Euler's pentagonal number theorem (\ref{euler2}).

\begin{thm}\label{euler5}
For any positive integer $n$, we have
$$
\sum\limits_{\lambda \in P_{do}(n)}(-1)^{\ell(\lambda)}=\left\{
\begin{array}{ll}
(-1)^k,  & \mbox{if}\ n=k^2,\\[6pt]
0, & \mbox{otherwise}.
\end{array}
\right.
$$
\end{thm}

Let $R_o(n)$ (resp. $R_e(n)$) denote the number
of partitions of $n$ into distinct parts such that the smallest part is odd,
and the number of parts is odd (resp. even). Then Theorem \ref{euler5} can
be restated as a theorem of Alladi \cite{alld092}.

\begin{thm}\label{add5}
For any positive integer $n$, we have
\begin{equation}\label{add6}
R_e(n)-R_o(n)=\left\{
\begin{array}{ll}
(-1)^k,  & \mbox{if}\ n=k^2,\\[6pt]
0, & \mbox{otherwise}.
\end{array}
\right.
\end{equation}
\end{thm}

Clearly, the problem of Pak is equivalent to finding a combinatorial proof of the
above theorem.
For example,
$n=9$, the $2$ partitions counted by $R_e(9)$ are
$$8+1, \quad 6+3,$$
while the $3$ partitions enumerated by $R_o(9)$ are
$$9, \quad
6+2+1, \quad 5+3+1.$$
Under the involution $\Psi$, the partitions are paired as follows
$$8+1 \leftrightarrow 9, \quad  6+3
\leftrightarrow 6+2+1,$$ while the triangular partition $5+3+1$ remains fixed.

\section{Alladi's partition theorem}

In this section, we apply the involution $\Psi$ presented in the previous section
to give a combinatorial interpretation of  a weighted partition theorem
derived by Alladi \cite{alld091} from Ramanujan's partial
theta identity (\ref{iden1}). While Berndt, Kim and Yee \cite{berndt09} constructed an
involution for the identity (\ref{iden1}),  they raised the question of how to
    translate their involution into a combinatorial proof of Alladi's weighted partition theorem. Even though our involution is not a direct answer to their question,
    it is likely that there is no easy way to make the translation. If so, our
    combinatorial interpretation can be considered as an indirect answer to the question of Berndt, Kim and Yee. The theorem of Alladi is stated as follows.

\begin{thm}\label{weig1}
For $\lambda=(\lambda_1,\lambda_2,\ldots,\lambda_l)\in P_{do}$, define
$\delta_i$ to be the least integer $\geq (\lambda_i-\lambda_{i+1})/2$, where
$\lambda_{l+1}$ is defined to be 0. Define the weight of $\lambda$ by
\begin{equation}\label{alldi1}
\omega_g(\lambda)=(-1)^l\prod\limits_{i=1}^la^{\delta_i}.
\end{equation}
Then we have
\begin{equation}\label{alldi2}
\sum\limits_{\lambda \in P_{do}(n)}\omega_g(\lambda)=\left\{
\begin{array}{ll}
(-a)^k,  & \mbox{if}\ n=k^2,\\[6pt]
0, & \mbox{otherwise}.
\end{array}
\right.
\end{equation}
\end{thm}

\pf For $\lambda \in P_{do}(n)$,
let $(\pi, \sigma)$ be the pair of partitions obtained from $\lambda$
 in Step $1$ of the Franklin type involution $\Psi$. It can be
 seen that
the exponent of $a$ in $\omega_g(\lambda)$ equals $\ell(\pi)+\ell(\lambda)$.
It is also clear that the quantity $\ell(\pi)+\ell(\lambda)$ remains unchanged in
Step $2$, that is
$$\ell(\pi)+\ell(\lambda)=\ell(\pi^*)+\ell(\lambda^*).$$
Let $\lambda^{*}=\pi^*+c_2(\sigma^*)$ in Step $3$, then the exponent of $a$ in $\omega_g(\lambda^*)$ equals $\ell(\pi^*)+\ell(\lambda^*)$.
Thus the involution $\Psi $ preserves the exponent of $a$ in $\omega_g(\lambda)$.
In view of the property (\ref{add8}), we see that $\omega_g(\lambda)$
and $\omega_g(\lambda^*)$ have  opposite signs.
Therefore, the partitions $\lambda$ in $P_{do}(n)$ cancel each except for the
the partition $\lambda=(2k-1,2k-3,\ldots,3,1)$ which has weight $(-a)^k$ for $n=k^2$.
This completes the proof. \qed

For example,
when $n=10$, there are six partitions in $P_{do}(10)$:
$$9+1, \quad 7+3, \quad 4+3+2+1, $$
$$7+2+1, \quad 6+3+1, \quad 5+4+1.$$
The involution $\Psi$ gives the following correspondence
$$9+1 \leftrightarrow 7+2+1, \quad  7+3
\leftrightarrow 6+3+1, \quad 4+3+2+1
\leftrightarrow 5+4+1.$$
Meanwhile, the weights of the partitions are given by
$$\omega_g(9+1)=a^5,\quad \omega_g(7+3)=a^4, \quad \omega_g(4+3+2+1)=a^4,$$
and
$$\omega_g(7+2+1)=-a^5,\quad \omega_g(6+3+1)=-a^4, \quad \omega_g(5+4+1)=-a^4.$$

\section{Another partition theorem of Alladi}

As will be seen,  the Franklin type involution $\Psi$ can be used to give a
combinatorial interpretation of
another weight partition theorem of Alladi by adopting a different weight assignment.
Alladi \cite{alld092} found a simple $q$-theoretic proof of
Andrews' partial theta identity
\begin{equation}\label{iden2x}
\sum\limits_{n=0}^{\infty}q^{2n}(q^{2n+2};q^2)_{\infty}(aq^{2n+1};q^2)_{\infty}
=\sum\limits_{k=0}^{\infty}(-a)^kq^{k^2},
\end{equation}
which is (\ref{iden2}) as mentioned before.
For $\lambda \in P_{do}$,  using a  simpler weight function
\begin{equation}\label{omega2}
 \omega_o(\lambda)=(-1)^la^{\ell_o(\lambda)},
\end{equation}
from the identity (\ref{iden2x})  Alladi deduced the following weighted
partition theorem.

\begin{thm} \label{weig2} Assume that the weight of a partition $\lambda\in P_{do}(n)$ is given by (\ref{omega2}). Then we have
\begin{equation}\label{alldi4}
\sum\limits_{\lambda \in P_{do}(n)}\omega_o(\lambda)=\left\{
\begin{array}{ll}
(-a)^k,  & \mbox{if}\ n=k^2,\\[6pt]
0, & \mbox{otherwise}.
\end{array}
\right.
\end{equation}
\end{thm}
\pf Let $\lambda \in P_{do}(n)$.
By the property (\ref{add8}),
it is easily seen that
the involution
$\Psi$ changes the number of even parts of $\lambda$ by $1$ and preserves the number of odd parts.
Consequently, the involution $\Psi$  preserves the exponent of $a$ given in the weight  $\omega_o(\lambda)$
and reverses the sign of $\omega_o(\lambda)$.  When $n$ is a square, say, $n=k^2$,
there exists exactly one partition which is undefined for $\Psi$,
that is $\lambda=(2k-1,2k-3,\ldots,3,1)$ whose weight equals $(-a)^k$.
This completes the proof.  \qed

We note that Theorem \ref{weig2}  can be translated back to the following identity:
\begin{equation}\label{add1}
\sum\limits_{n=1}^{\infty}-aq^{2n-1}(q^{2n};q^2)_{\infty}(aq^{2n+1};q^2)_{\infty}=
\sum\limits_{n=1}^{\infty}(-a)^nq^{n^2},
\end{equation}
which takes a different form compared with the identity (\ref{iden2x}). Nevertheless, as shown by
Alladi \cite{alld092}, (\ref{add1}) can be deduced from (\ref{iden2x}).

\section{A partition theorem related to Andrews' identity}

As we have seen in the previous section,
Theorem \ref{weig2} is a direct translation of the identity (\ref{add1})
rather than the idenity (\ref{iden2}). So we are led to a partition identity
 directly derived from the identity (\ref{iden2}) which can be  proved
with the aid of our involution $\Psi$.
Recall that $Q$ denotes the set of partitions into distinct non-negative parts
with  the smallest part being even. Let $Q(n)$ denote such partitions of $n$ in $Q$.
 For a partition
$\lambda=(\lambda_1,\lambda_2,\ldots,\lambda_l)\in Q$,
define the weight of $\lambda$ by
\begin{equation} \label{omega3}
\omega_e(\lambda)=(-1)^{l-1}a^{\ell_o(\lambda)}.
\end{equation}
Then we have the following partition identity.

\begin{thm}\label{weig4} We have
\begin{equation}\label{andr2}
\sum\limits_{\lambda \in Q(n)}\omega_e(\lambda)=\left\{
\begin{array}{ll}
(-a)^k,  & \mbox{if}\ n=k^2,\\[6pt]
0, & \mbox{otherwise}.
\end{array}
\right.
\end{equation}
\end{thm}

\pf Let $\lambda$ be a partition in $Q(n)$.
Let $s(\lambda)$
denote the smallest part of the partition $\lambda$, and
let $ss(\lambda)$ denote the second small
part of $\lambda$. Define an involution  $\psi$  by the following procedure.
Three cases are considered.
\begin{itemize}
\item[(i)] Assume that $s(\lambda)=0$ and $ss(\lambda)$ is even. Delete the part $s(\lambda)$ in $\lambda$ and denote the resulting partition  by $\lambda^{*}$.
    It can be seen that $\lambda^{*}\in Q(n)$.
    \item[(ii)] Assume that  $s(\lambda) \neq
0$. Add $0$ to $\lambda$ as a new part. Denote the resulting
partition by $\lambda^{*}$. Then we have $\lambda^{*}\in Q(n)$.
\item[(iii)] Assume that
$s(\lambda)=0$ and $ss(\lambda)$ is odd.
In this case, $\lambda$  can be considered as a partition in $P_{do}(n)$ by disregarding the zero part $s(\lambda)$ so that we can apply $\Psi$ to $\lambda$.
\end{itemize}

According to the above construction, $\psi$ is  a sign-reversing and weight-preserving involution for which the partition
 $\lambda=(2k-1,2k-3,\ldots,3,1,0)\in Q(n)$ is defined as the fixed point for $n=k^2$. This
 completes the proof. \qed

For example,
$n=10$, there are fourteen partitions in $Q(10)$:
$$10, \quad 8+2, \quad 6+4, \quad 5+3+2$$
$$10+0, \quad 8+2+0, \quad 6+4+0, \quad 5+3+2+0$$
$$9+1+0, \quad 7+3+0, \quad 4+3+2+1+0, $$
$$7+2+1+0, \quad 6+3+1+0, \quad 5+4+1+0.$$
In this example, the involution  $\psi$ gives the following correspondence
$$10 \leftrightarrow 10+0, \quad  8+2
\leftrightarrow 8+2+0, \quad 6+4
\leftrightarrow 6+4+0, \quad 5+3+2
\leftrightarrow 5+3+2+0$$
$$9+1+0 \leftrightarrow 7+2+1+0, \,  7+3+0
\leftrightarrow 6+3+1+0, \, 4+3+2+1+0
\leftrightarrow 5+4+1+0.$$
The weights of partitions in $Q(10)$ are listed below, and it can be seen that
$\psi$ is indeed weight-preserving and sign-reversing,
$$\omega_e(10)=1,\quad \omega_e(8+2)=-1, \quad \omega_e(6+4)=-1, \quad \omega_e(5+3+2)=a$$
$$\omega_e(9+1+0)=a^5,\quad \omega_e(7+3+0)=a^4, \quad \omega_e(4+3+2+1+0)=a^4,$$
$$\omega_e(10+0)=-1,\quad \omega_e(8+2+0)=1, \quad \omega_e(6+4+0)=1, \quad \omega_e(5+3+2+0)=-a$$
$$\omega_e(7+2+1+0)=-a^5,\quad \omega_e(6+3+1+0)=-a^4, \quad \omega_e(5+4+1+0)=-a^4.$$

\section{Connection to an identity of Andrews}

In this section, we consider the special case of Theorem \ref{weig4} when setting
$a=-1$, that is,
\begin{equation}\label{iden4}
\sum\limits_{n=0}^{\infty}q^{2n}(q^{2n+2};q^2)_{\infty}(-q^{2n+1};q^2)_{\infty}
=\sum\limits_{k=0}^{\infty}q^{k^2},
\end{equation}
 which turns out to be related to
a problem  proposed by Andrews \cite{andr72} in 1972, see also, Andrews \cite[pp. 156-157]{andr76}. The original problem of Andrews is stated below.

{\noindent \bf  A Problem of Andrews.} Let $q_e(n)$ (resp.
$q_o(n)$) denote the number of partitions in $Q(n)$ that have an even number (resp. odd number) of even parts. Prove that
\begin{equation}\label{iden5}
q_o(n)-q_e(n)=\left\{
\begin{array}{ll}
1,  & \mbox{if}\ n=k^2,\\[6pt]
0, & \mbox{otherwise}.
\end{array}
\right.
\end{equation}

Clearly, the left hand side of (\ref{iden4}) counts the number of
partitions $\lambda$ in $Q$ with the sign $(-1)^{\ell_e(\lambda)-1}$
attached to $\lambda$. The sign
$(-1)^{\ell_e(\lambda)-1}$ equals the weight of $\lambda$ by setting $a=-1$ in (\ref{omega3}),
namely,
\[ \omega_e(\lambda)=(-1)^{l-1}a^{\ell_o(\lambda)}. \]
Thus we can apply the above involution $\psi$
 defined in the previous section give a combinatorial
interpretation of the identity (\ref{iden5}).

When $a=-1$, the identity (\ref{andr2}) can be rewritten as
\begin{equation}\label{iden11}
\begin{split}
&\sum\limits_{\lambda \in Q(n)}\omega_e(\lambda)=
\sum\limits_{\lambda \in Q(n)}(-1)^{l-1}(-1)^{\ell_o(\lambda)}=\sum\limits_{\lambda \in Q(n)}(-1)^{\ell_e(\lambda)-1}\\
&\quad = q_o(n)-q_e(n) =\left\{
\begin{array}{ll}
1,  & \mbox{if}\ n=k^2,\\[6pt]
0, & \mbox{otherwise}.
\end{array}
\right.
\end{split}
\end{equation}
It is clear from  (\ref{add8}) that the
involution $\psi$ only changes the number of even parts of $\lambda \in Q(n)$ by
$1$.
Thus  the identity (\ref{iden5}) follows from the involution $\psi$.

For example, when $n=9$, the five partitions enumerated by $q_e(9)$
are
\[ 8+1+0,\quad  7+2+0, \quad 6+3+0,\quad 5+4+0, \quad 4+3+2,\]
and the six
partitions enumerated by $q_o(9)$ are
\[ 9+0, \quad 7+2, \quad 6+2+1+0, \quad
5+4, \quad 5+3+1+0, \quad 4+3+2+0.\]
Under the involution $\psi$, the partitions are paired as follows
\[ 8+1+0 \leftrightarrow 9+0,\quad  7+2+0 \leftrightarrow 7+2,\quad  6+3+0
\leftrightarrow 6+2+1+0,   \]
\[  5+4+0 \leftrightarrow 5+4,\quad  4+3+2
\leftrightarrow 4+3+2+0.\]
The partition $5+3+1+0$ is the fixed point.

To conclude this section, we  remark that although the above identity (\ref{iden5})
takes a slightly different form with
Theorem \ref{add5}, it is obvious that they are equivalent, as noted by Alladi \cite{alld091}.
The involution $\psi$ gives an explanation of the
 equivalence these two identities.

\section{A more general partition theorem}

In this section, we present the following weighted
form of Andrews' identity (\ref{add7}):
\begin{equation}\label{add9}
\begin{split}
&\sum\limits_{n=0}^{\infty}q^{2mn}(q^{2mn+2m};q^{2m})_{\infty}(aq^{2mn+1};q^2)_{\infty}\\
&=1+\sum\limits_{k=1}^{\infty}(-a)^kq^{k^2}\prod\limits_{j=1}^k(1+q^{2j}+q^{4j}+\cdots+q^{2(m-1)j}),
\end{split}
\end{equation}
which reduces to (\ref{add7}) by setting $a=-1$ and
reduces to (\ref{iden2}) by setting $m=1$.
By extending our Franklin type involution $\Psi$, we obtain a combinatorial
interpretation of the above generalization.

Let us introduce some notation.
For a positive integer $m$, let $A_{k,\,m}$ denote the set of
partitions into $k$ distinct nonnegative parts such that all the
even parts are multiples of $2m$ and the smallest part is even. Let
$A_m=\bigcup \limits_{k=0}^{\infty}A_{k,\,m}$ and let $A_m(n)$ be the set of
partitions of $n$ in $A_m$. In this notation,
the generating function for partitions $\lambda \in A_m$ equals
\begin{equation}\label{iden12}
\sum\limits_{n=0}^{\infty}q^{2mn}(-q^{2mn+2m};q^{2m})_{\infty}(-q^{2mn+1};q^2)_{\infty}.
\end{equation}
To give a combinatorial interpretation of the right hand side of
(\ref{add9}), let $H_{k,\,m}$ denote the set of partitions $\lambda_{k,\,m}$
such that each
part of $\lambda_{k,\,m}$ is less than or equal to $k$ and the
multiplicity of each part is an even number less than $2m$.
Then the generating function for  partitions
$\lambda_{k,\,m}$ in $H_{k,\,m}$ equals
$$\prod\limits_{j=1}^k(1+q^{2j}+q^{4j}+\cdots+q^{2(m-1)j}).$$
The factor $q^{k^2}$ equals the generating function of the triangular partition
 \[ T_k=
 (2k-1,2k-3,\ldots,3,1).\]

In order to give a combinatorial explanation of the identity (\ref{add9}),
we shall give another interpretation of the right hand side of (\ref{add9}).
To this end, let $B_{k,\,m}$ denote the set of partitions
$\pi=(\pi_1,\pi_2,\ldots,\pi_k)$ into distinct odd parts such that
the difference of consecutive parts is less than or equal to $2m$,
namely, $\pi_i-\pi_{i+1}\leq 2m$ for $1\leq i\leq k$ with the
convention that $\pi_{k+1}=0$. Set \[ B_m=\bigcup\limits_{k=0}^{\infty}B_{k,\,m},\]
and let $B_m(n)$ be the set of
partitions of $n$ in $B_m$.
Then we have the following correspondence.

\begin{thm}\label{thm8}
There exists a bijection between the set $B_{k,\,m}$ and the set $\{T_k\} \times H_{k,\,m}$.
\end{thm}
\pf We proceed to construct a  bijection from $B_{k,\,m}$ to
$\{T_k\} \times H_{k,\,m}$.
For a partition $\lambda=(\lambda_1,\lambda_2,\ldots,\lambda_k)\in
B_{k,\,m}$, using  the  Ferrers diagram, we can generate the triangular partition $T_k$
and a partition $\lambda_{k,\,m}\in H_{k,\,m}$ by the following
procedure. Let $i$ be the
largest integer such that $\lambda_i-\lambda_{i+1}=2j>2$ and $j\leq m$ with the
convention that $\lambda_{k+1}=0$. Then we remove $2(j-1)$ columns
of length $i$ from $\lambda$ and add them to $\lambda_{k,\,m}$ as rows.
Repeating this procedure until there does not exist such $i$.
Finally, the remaining partition is the triangular partition $T_k$.
 It can be seen that
$\lambda_{k,\,m}\in  H_{k,\,m}$.

This process is reversible. Given the triangular $T_k$ and a partition $ \lambda_{k,\,m}\in H_{k,\,m}$,
let $\lambda=T_k+\lambda^{'}_{k,\,m}$. Then, we have $\lambda \in
B_{k,\,m}$. This completes the proof. \qed

Below is  an example when $\lambda=(19,15,9,5,3)\in B_{5,\,3}$.

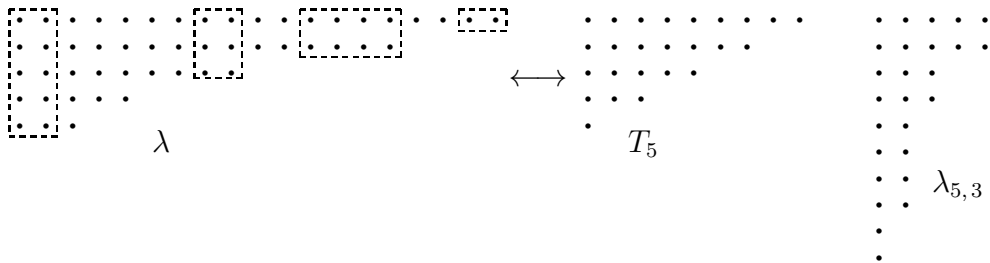
\begin{center}
\begin{picture}(375,90)
\multiput(0,90)(10,0){19}{\circle*{2}}\multiput(0,80)(10,0){15}{\circle*{2}}
\multiput(0,70)(10,0){9}{\circle*{2}}\multiput(0,60)(10,0){5}{\circle*{2}}
\multiput(0,50)(10,0){3}{\circle*{2}}
\put(185,65){$\longleftrightarrow$}
\multiput(215,90)(10,0){9}{\circle*{2}}\multiput(215,80)(10,0){7}{\circle*{2}}
\multiput(215,70)(10,0){5}{\circle*{2}}\multiput(215,60)(10,0){3}{\circle*{2}}
\multiput(215,50)(10,0){1}{\circle*{2}}

\multiput(325,90)(10,0){5}{\circle*{2}}\multiput(325,80)(10,0){5}{\circle*{2}}
\multiput(325,70)(10,0){3}{\circle*{2}}\multiput(325,60)(10,0){3}{\circle*{2}}
\multiput(325,50)(10,0){2}{\circle*{2}}\multiput(325,40)(10,0){2}{\circle*{2}}
\multiput(325,30)(10,0){2}{\circle*{2}}\multiput(325,20)(10,0){2}{\circle*{2}}
\put(325,10){\circle*{2}}\put(325,0){\circle*{2}}

\multiput(-4,94)(0,-4){12}{\line(0,-1){2}}\multiput(14,94)(0,-4){12}{\line(0,-1){2}}
\multiput(-4,94)(4,0){5}{\line(1,0){2}}\multiput(-4,46)(4,0){5}{\line(1,0){2}}

\multiput(66,94)(0,-4){7}{\line(0,-1){2}}\multiput(84,94)(0,-4){7}{\line(0,-1){2}}
\multiput(66,94)(4,0){5}{\line(1,0){2}}\multiput(66,68)(4,0){5}{\line(1,0){2}}

\multiput(106,94)(0,-4){5}{\line(0,-1){2}}\multiput(144,94)(0,-4){5}{\line(0,-1){2}}
\multiput(106,94)(4,0){10}{\line(1,0){2}}\multiput(106,76)(4,0){10}{\line(1,0){2}}

\multiput(166,94)(0,-4){2}{\line(0,-1){2}}\multiput(184,94)(0,-4){2}{\line(0,-1){2}}
\multiput(166,94)(4,0){5}{\line(1,0){2}}\multiput(166,86)(4,0){5}{\line(1,0){2}}

\put(50,40){$\lambda$}\put(230,40){$T_5$}\put(345,25){$\lambda_{5,\,3}$}
\end{picture}
\end{center}

From Theorem
\ref{thm8}, we conclude that the generating function for partitions
$\lambda \in B_m$ equals
\begin{equation}\label{iden13}
1+\sum\limits_{k=1}^{\infty}q^{k^2}
\prod\limits_{j=1}^k(1+q^{2j}+q^{4j}+\cdots+q^{2(m-1)j}).
\end{equation}
Using the identities (\ref{iden12}) and (\ref{iden13}),
we obtain the the following number-theoretic interpretation of the identity
(\ref{add9}) in terms of weighted partitions.

\begin{thm}\label{weig5}
Assume that $m\geq 1$. For
$\lambda=(\lambda_1,\lambda_2,\ldots,\lambda_l)\in A_m$, define the weight of
$\lambda$ by
\begin{equation}\label{andr3}
\omega_1(\lambda)=(-1)^{l-1}a^{\ell_o(\lambda)}.
\end{equation}
On the other hand, for
$\mu=(\mu_1,\mu_2,\ldots,\mu_l)\in B_m$, define the weight of $\mu$ by
\begin{equation}\label{andr4}
\omega_2(\mu)=(-a)^l.
\end{equation}
Then the following relation holds
\begin{equation}\label{andr5}
\sum\limits_{\lambda \in A_m(n)}\omega_1(\lambda)=
\sum\limits_{
\mu \in B_m(n)}\omega_2(\mu).
\end{equation}
\end{thm}

Since $A_1=Q$ and $B_1$ consists of only triangular partitions  $T_k=(2k-1,2k-3,\ldots,3,1)$,
 Theorem \ref{weig5} reduces to Theorem \ref{weig4}   by  setting
$m=1$.

The proof of Theorem \ref{weig5} relies on the notion of $2m$-modular diagrams,
see  \cite{pak06}. Recall that  the $2m$-modular diagram
of a partition $\lambda$ is
defined to be Young diagram by placing the integer $2m$  in the squares of each row
possibly except for the last square, and the last square of each row may be filled
with an integer not exceeding $2m$.

Let $P_{do}^m(n)$ denote the set of partitions of $n$ into distinct parts
such that all the even parts are multiples of $2m$ and the smallest
part is odd. Using the $2m$-modular diagrams of partitions, we can  extend
the Franklin type involution $\Psi$ on $P_{do}(n)$ to
 $P_{do}^m(n)$, and we denote it by $\Psi_m$. The explicit construction of
 $\Psi_m$ is analogous to the three steps of the involution
 $\Psi$ in Section 3, and hence it is omitted.
 Furthermore, we can extend the involution $\psi$ on $Q(n)$
 to $A_m(n)$ with the aid of
$\Psi_m$ to give a combinatorial proof of Theorem \ref{weig5}.
 Since the proof of Theorem \ref{weig5} is similar to that of Theorem \ref{weig4},
it is also omitted. Here is an example of the involution $\psi_m$ for $m=2$.
For $\lambda=(20,16,11,5,3,0)\in A_2(55)$, we have
$\psi_2(\lambda)=(20,19,13,3,0)\in A_2(55)$. The following figure is an
illustration of the procedure to construct  $\psi_2(\lambda)$.

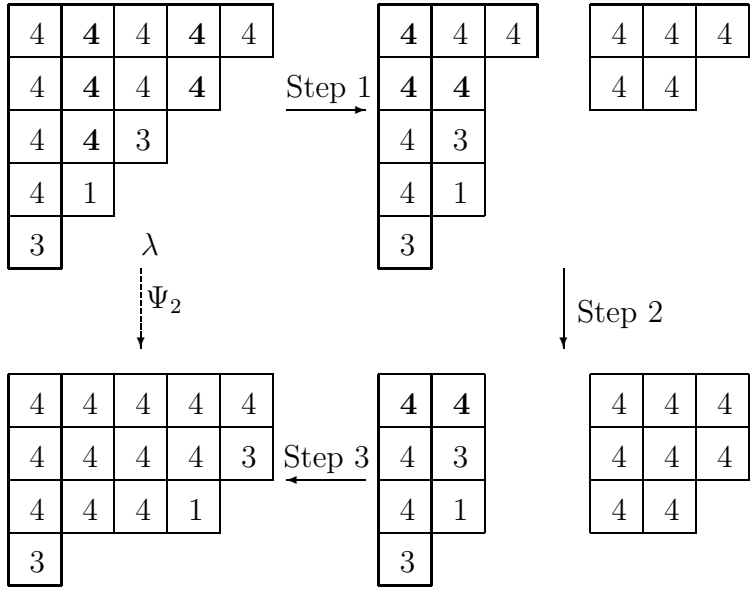
\begin{center}
\begin{picture}(280,240)
\put(0,240){\line(1,0){100}}\put(0,220){\line(1,0){100}}
\put(0,200){\line(1,0){80}}\put(0,180){\line(1,0){60}}
\put(0,160){\line(1,0){40}}\put(0,140){\line(1,0){20}}
\put(0,240){\line(0,-1){100}}\put(20,240){\line(0,-1){100}}
\put(40,240){\line(0,-1){80}}\put(60,240){\line(0,-1){60}}
\put(80,240){\line(0,-1){40}}\put(100,240){\line(0,-1){20}}
\multiput(8,225)(0,-20){4}{4}\multiput(28,225)(0,-20){3}{\bf 4}
\multiput(48,225)(0,-20){2}{4}\multiput(68,225)(0,-20){2}{\bf 4}
\put(88,225){4}\put(48,185){3}\put(28,165){1}\put(8,145){3}

\put(105,200){\vector(1,0){30}}\put(105,205){Step $1$}

\put(140,240){\line(1,0){60}}\put(140,220){\line(1,0){60}}
\put(140,200){\line(1,0){40}}\put(140,180){\line(1,0){40}}
\put(140,160){\line(1,0){40}}\put(140,140){\line(1,0){20}}
\put(140,240){\line(0,-1){100}}\put(160,240){\line(0,-1){100}}
\put(180,240){\line(0,-1){80}}\put(200,240){\line(0,-1){20}}
\put(148,225){\bf 4}\multiput(168,225)(20,0){2}{4}\multiput(148,205)(20,0){2}{\bf 4}
\multiput(148,185)(0,-20){2}{4}
\put(168,185){3}\put(168,165){1}\put(148,145){3}

\put(220,240){\line(1,0){60}}\put(220,220){\line(1,0){60}}
\put(220,200){\line(1,0){40}}
\put(220,240){\line(0,-1){40}}\put(240,240){\line(0,-1){40}}
\put(260,240){\line(0,-1){40}}\put(280,240){\line(0,-1){20}}
\multiput(228,225)(20,0){3}{4}\multiput(228,205)(20,0){2}{4}

\put(210,140){\vector(0,-1){30}}\put(215,120){Step $2$}

\put(140,100){\line(1,0){40}}\put(140,80){\line(1,0){40}}
\put(140,60){\line(1,0){40}}\put(140,40){\line(1,0){40}}
\put(140,20){\line(1,0){20}}
\put(140,100){\line(0,-1){80}}\put(160,100){\line(0,-1){80}}
\put(180,100){\line(0,-1){60}}
\multiput(148,85)(20,0){2}{{\bf 4}}
\multiput(148,65)(0,-20){2}{4}
\put(168,65){3}\put(168,45){1}\put(148,25){3}

\put(220,100){\line(1,0){60}}\put(220,80){\line(1,0){60}}
\put(220,60){\line(1,0){60}}\put(220,40){\line(1,0){40}}
\put(220,100){\line(0,-1){60}}\put(240,100){\line(0,-1){60}}
\put(260,100){\line(0,-1){60}}\put(280,100){\line(0,-1){40}}
\multiput(228,85)(20,0){3}{4}\multiput(228,65)(20,0){3}{4}
\multiput(228,45)(20,0){2}{4}

\put(135,60){\vector(-1,0){30}}\put(100,65){ Step $3$}

\put(0,100){\line(1,0){100}}\put(0,80){\line(1,0){100}}
\put(0,60){\line(1,0){100}}\put(0,40){\line(1,0){80}}
\put(0,20){\line(1,0){20}}
\put(0,100){\line(0,-1){80}}\put(20,100){\line(0,-1){80}}
\put(40,100){\line(0,-1){60}}\put(60,100){\line(0,-1){60}}
\put(80,100){\line(0,-1){60}}\put(100,100){\line(0,-1){40}}
\multiput(8,85)(20,0){5}{4}\multiput(8,65)(20,0){4}{4}
\multiput(8,45)(20,0){3}{4}
\put(88,65){3}\put(68,45){1}\put(8,25){3}

\put(50,145){$\lambda$}\multiput(50,140)(0,-3){10}{\line(0,-1){2}}
\put(50,115){\vector(0,-1){5}}\put(52,125){$\Psi_2$}

\end{picture}
\end{center}

\vskip 3mm \noindent {\bf Acknowledgments.} This work was
supported by the 973 Project, the PCSIRT Project of the Ministry
of Education, and  the National
Science Foundation of China.

\end{document}